\documentclass[11pt]{amsart}
\usepackage{amsfonts}
\usepackage{amssymb}
\usepackage[arrow,matrix]{xy}
\usepackage{amsmath,amssymb, bbm, amscd, amsthm, mathrsfs, hyperref, mathtools}
\usepackage[titletoc]{appendix}
\usepackage{tikz-cd}

\DeclareMathAlphabet{\mathpzc}{OT1}{pzc}{m}{it}

\theoremstyle{plain}
\newtheorem{thm}{Theorem}[section]

\newtheorem{lem}[thm]{Lemma}
\newtheorem{prop}[thm]{Proposition}

\theoremstyle{definition}

\newtheorem{eg}[thm]{Example}

\newtheorem{rem}[thm]{Remark}
\newtheorem{con}[thm]{Convention}

\numberwithin{equation}{section}

\def\mc{\mathcal}

\def\it{\textit}

\def\ot{\otimes}

\def\ra{\rightarrow}

\def\id{\operatorname {id}}

\def\End{\operatorname {End}}

\def\Sym{\operatorname {Sym}}

\def\kk{\mathbbm{k}}

\def\NN{\mathbb{N}}

\begin{document}
	\title[Isotopes and Zhang twists]{\bf  Isotopes of biracks and Zhang twists of  algebras}

	\author{Xiaolan YU}
	\address {Xiaolan YU\newline Department of Mathematics, Hangzhou Normal University, Hangzhou, Zhejiang 310036, China}
	
	\email{xlyu@hznu.edu.cn}

\author{Yanfei ZHANG}
\address {Yanfei ZHANG\newline Department of Mathematics, Hangzhou Normal University, Hangzhou, Zhejiang 310036, China}
	
	\email{zyfeig@163.com}
	
	\date{}

	\begin{abstract} In this paper,  we introduce the notion of an  $\NN^p$-graded birack and construct its isotope. Every involutive $\NN^p$-graded birack gives rise to an $\NN^p$-graded Yang-Baxter algebra. We study the relation between isotopes of   involutive $\NN^p$-graded biracks and  Zhang twists of $\NN^p$-graded Yang-Baxter algebras. As an example,  Yang–Baxter algebras determined by distributive solutions  are proved to be Zhang twists of  polynomial algebras.

	\end{abstract}

	\maketitle
	
	\section*{Introduction}\label{0}
	The  Yang–Baxter equation is a cornerstone of mathematical physics that arises
	in the study of integrable systems, quantum groups, and quantum information theory.   Its origins trace back to foundational papers by Yang \cite{yang} and Baxter \cite{bax}. Let $V$ be a vector space. A solution
	of the Yang–Baxter equation is a linear mapping $R : V \ot V \ra V \ot V$ 
	such that the equality	
$$(\id \otimes R)(R \otimes \id)(\id \otimes R)=(R \otimes \id)(\id \otimes R)(R \otimes \id)$$
		of operators hold in the space  $\End(V \otimes V \otimes V)$. The problem of describing all possible solutions is highly complex and therefore there
		were some simplifications introduced.	
		
		A particularly nice class of solutions is provided by set-theoretic solutions proposed by Drinfeld \cite{dr}. Let $X$ be a set and $r:X\times X\ra X\times X$ a bijective mapping. $r$ is called a set-theoretic solution of the Yang-Baxter equation if the following braid relation holds
		$$(\id \times r)(r \times \id)(\id \times r)=(r \times \id)(\id \times r)(r \times \id).$$
			It is clear that a set-theoretic solution extends to a linear one, but more important than this is
	that set-theoretic solutions lead to their own remarkable algebraic and combinatoric structure. During the last three decades, the study of set-theoretic solutions and related structures has notably intensified, see for instance \cite{do}, \cite{dos}, \cite{ess}, \cite{gi}, \cite{gi3}, \cite{gim1}, \cite{gim2}, \cite{giv}, \cite{guv}, \cite{lyz}.

	On the algebra side, each set-theoretic solution $(X,r)$ associates a quadratic Yang-Baxter algebra $A(X,r)$ as proposed in \cite[Section 6]{gim2}. If $X$ is finite and $r$ is nondegenerate and involutive, then  $A(X,r)$ has remarkable algebraic, homological and combinatorial properties \cite{gi2}, \cite{giv}. 	Our motivation of this paper is to try to  find more ways to study  properties of Yang-Baxter algebras.  
	
	It is known that there is a one-to-one correspondence between non-degenerate, involutive set-theoretic solutions of the Yang-Baxter equation and involutive biracks. A birack $(X, \circ, \backslash_\circ, \bullet, /_\bullet)$ is an  algebra such that $(X, \circ, \backslash_\circ)$ is a left quasigroup, $(X, \bullet, /_\bullet)$ is a right quasigroup, and some additional identities are satisfied. Such correspondence allows us to study solutions of the Yang-Baxter equations via involutive biracks.	
	
	In quasigroup theory (see, e.g., \cite[Section II.2]{pf}), isotopy is a standard method for transforming  one quasigroup into another. In their work of  involutive set-theoretic solutions of the Yang-Baxter equation of multipermutation level 2,  Jedli\v{c}ka, Pilitowska, and Zamojska-Dzienio used a special form of isotope  to construct non-distributive solutions from distributive solutions \cite{jpz}. Inspired by their work, in Section \ref{S2}, we introduce the notions of an $\NN^p$-graded birack and its  isotope, and then discuss their properties. Our main theorem is presented in Section \ref{S3}. Each involutive  $\NN^p$-graded birack $(X, \circ, \backslash_\circ, \bullet, /_\bullet)$  determines an  $\NN^p$-graded Yang-Baxter algebra $A(X,{\bf R}^\circ )$. We find that under some conditions, the algebra  $A^*(X,{\bf R}^* )$ determined by an isotope  $(X, *, \backslash_*, \diamond, /_\diamond)$ of the birack $(X, \circ, \backslash_\circ, \bullet, /_\bullet)$ is isomorphic to a Zhang twist of the algebra $A(X,{\bf R}^\circ )$ (Theorem \ref{T1}). The Zhang twist is a twist of graded algebras introduced by J.J. Zhang in order to describe graded algebras whose graded module categories are equivalent \cite{zh1}. Many important properties of graded algebras, including Gelfand–Kirillov dimension, the Noetherian property, being a domain, global dimension, and Artin–Schelter regularity are preserved under Zhang twist.
As an example our main theorem, we show that the Yang–Baxter algebras arising from distributive solutions are Zhang twists of polynomial algebras. This class of algebras is believed to share properties similar to those of polynomial algebras. We will discuss them in detail in a subsequent paper.

	\section{Preliminaries}\label{S1}

	Throughout $\kk$ is a base field of characteristic zero. All vector spaces and algebras are over $\kk$.

	\subsection{Quasigroups and biracks}
	A \it{left quasigroup} is an algebra $(X, \circ, \backslash_\circ)$ with two binary operations: $\circ$ and 
	$\backslash_\circ$, such that for every $x, y \in X $, the
	following conditions holds
	\begin{equation}\label{lq}
		x \circ(x \backslash _\circ y)=y=x \backslash_ \circ(x \circ y).
	\end{equation}

	A \it{right quasigroup} is defined analogously as an algebra $(  X, \bullet, /_ \bullet  )$ with two binary operations: $\bullet$ and $/_ \bullet$, such that for every   $x, y \in X$,  the following conditions holds
	\begin{equation}\label{rq}
		(y /_\bullet x) \bullet x=y=(y \bullet x) /_\bullet x.
	\end{equation}

	Let $(X, \circ, \backslash_\circ)$ be a left quasigroup. For each $x\in X$, the \it{left translation}  $\mc{L}_x:X\ra X$ by $x$ is defined as 
	$$\mc{L}_x(a)=x\circ a,$$ for any $a\in X$. Condition (\ref{lq}) guarantees that each $\mc{L}_x$ is bijective, with $\mc{L}_x^{-1}(a)=x\backslash_\circ a$. Similarly, Condition 
	(\ref{rq}) gives that each \it{right translation} $\mc{R}_x:X\ra X$ by $x$ defined by $\mc{R}_x(a)=a\bullet x$ is bijective with $\mc{R}_x^{-1}(a)=a/_\bullet x$.
	
	For a left quasigroup $(X, \circ, \backslash_\circ)$, a bijection  $f:X\ra X$ is called an \it{automorphism} of $(X, \circ, \backslash_\circ)$ if
	$$f(x\circ y )=f(x)\circ f(y ),$$
	for every $x,y\in X$.	It is easy to see that a bijection  $f:X\ra X$ is an automorphism of $(X, \circ, \backslash_\circ)$
	if and only if  $$f\mc{L}_x f^{-1}=\mc{L}_{f(x)},$$ for every $x\in X$.
	
	An automorphism of a  right quasigroup can be defined analogously. A bijection  $f:X\ra X$ is an automorphism of $(X, \bullet, /_\bullet)$
	if and only if  $$f\mc{R}_x f^{-1}=\mc{R}_{f(x)},$$ for every $x\in X$.
	
	A left quasigroup $(X, \circ, \backslash_\circ)$ is called
	\begin{itemize}
		\item \it{non-degenerate} if the mapping
		$$\begin{array}{rcl}
			X&\longrightarrow& X\\
			x&\longmapsto&x\backslash_\circ x
		\end{array}
		$$	is a bijection;
		\item \it{right cyclic} if for every $x, y, z \in X$:
		$$\left(x \backslash_{\circ} y\right) \backslash_{\circ}\left(x \backslash_{\circ} z\right)=\left(y \backslash_{\circ} x\right) \backslash_{\circ}\left(y \backslash_{\circ} z\right),$$
		or equivalently 
		$$ \mc{L}_{x} \mc{L}_{x \backslash_{\circ} y}=\mc{L}_{y} \mc{L}_{y \backslash_{\circ} x}.$$
	\end{itemize}

	An algebra $(  X, \circ,\backslash_{\circ}, \bullet, /_\bullet  )$ with four binary operations is called a \it{birack}, if $(  X, \circ,\backslash_{\circ}  )$ is a left quasigroup, $(  X, \bullet, /_ \bullet  )$ is a right quasigroup and the following holds for any  $x, y, z \in X$,
	$$	\begin{aligned}
		x \circ(y \circ z) & =(x \circ y) \circ((x \bullet y) \circ z), \\
		(x \circ y) \bullet((x \bullet y) \circ z) & =(x \bullet(y \circ z)) \circ(y \bullet z), \\
		(x \bullet y) \bullet z & =(x \bullet(y \circ z)) \bullet(y \bullet z).
	\end{aligned}$$ 
	If the set $X$ is finite, we will call the birack \it{finite}. 
	
	\begin{con}
		All biracks  we study in this paper are finite.
	\end{con}

	A birack is called \it{involutive} if it additionally satisfies
	$$\begin{array}{l}
		
		(x \circ y) \circ(x \bullet y)=x \\
		(x \circ y) \bullet(x \bullet y)=y
	\end{array}$$
for every $x, y\in X$.

The following one to one correspondence  between involutive biracks and non-degenerate right cyclic left quasigroups  is well-known (see e.g. \cite[Proposition 1.5]{de}). 
	
	\begin{lem}\label{L4} \begin{enumerate}
			\item Let $\left(X, \circ, \backslash_{\circ}, \bullet, /_\bullet\right)$  be an involutive birack. Then $(X, \circ, \backslash_{\circ})$ is a non-degenerate right cyclic left quasigoup.
			\item Conversely, let $(X, \circ, \backslash_{\circ})$ be a non-degenerate right cyclic left quasigoup. For every  $x, y \in X$, define 
			$$x \bullet y=\mc{R}_{y}(x)=(x \circ y) \backslash_{\circ} x, $$
			$$x /_\bullet y=\mc{R}_{y}^{-1}(x).$$
			Then the algebra  $\left(X, \circ, \backslash_{\circ}, \bullet, / _\bullet\right)$  is an involutive birack. 
		\end{enumerate}
		
	\end{lem}
	
	The operations  $\backslash_{\circ}$  and $/_\bullet$ are connected by
	$$\left(x \backslash_{\circ} x\right) /_\bullet\left(x \backslash_{\circ} x\right)=x \quad \text { and } \quad(x /_\bullet x) \backslash_{\circ}(x /_\bullet x)=x.$$
	
	A birack $(  X, \circ,\backslash_{\circ}, \bullet, /_\bullet  )$ is said to satisfy the condition {\bf lri} provided that 
	$$(x \circ y) \bullet x=y=x \circ(y \bullet x),$$
	all $x,y\in X$. That is,   $\mc{R}_{x}=\mc{L}_{x}^{-1}$ for all $x\in X$.
	
	\begin{rem}
		An  involutive birack  $(X, \circ,\backslash_{\circ}, \bullet, /_\bullet  )$ satisfies Condition {\bf lri} if and only if 
		\begin{equation}
			(x \circ y) \backslash_{\circ} x=y\backslash_\circ x,
		\end{equation}
		for any $x,y\in X$.
	In fact, in an involutive birack  $x\bullet y=(x \circ y) \backslash_{\circ} x$. Hence, 
	$$
	\mc{R}_{y}(x)=\mc{L}_{y}^{-1}(x) \iff (x \circ y) \backslash_{\circ} x=y\backslash_\circ x.
	$$
		\end{rem}

	\subsection{Set-theoretic solution of of the Yang-Baxter equation} By a \it{quadratic set}, we mean a pair $(X,r)$, where $X$ is a nonempty set and  $r:X\times X\ra X\times X$  a bijective mapping. We write the image of $(x,y)$ under $r$ as 
	$$r(x, y)=\left({}^xy, x^y\right).$$
	The above equation defines a ``left action'' $\mc{L}:X\times X\ra X$, and a  ``right action'' $\mc{R}:X\times X\ra X$, on $X$ as: 
	$$\mc{L}_x(y)={}^xy,\;\;\;\;\;  \mc{R}_y(x)=x^y,$$
	for all $x,y\in X$.
	
	A quadratic set $(X,r)$ is called: 
	\begin{itemize}
		\item 	\it{non-degenerate} if the mappings $\mc{L}_x$ and $\mc{R}_x$ are  bijections, for all  $x\in X$;
		\item \it{involutive} if  $r^{2}={\id}_{X^{2}}$;
		\item \it{square free} if $ r(x, x)=(x, x)$, for every  $x \in X$;
		\item a \it{quantum binomial set} if it is nondegenerate, involutive and square-free.
	\end{itemize}
	
	A \it{set-theoretic solution of the Yang–Baxter equation} (YBE) is a quadratic set $(X, r)$, such that the following braid relation holds
	$$(\id \times r)(r \times \id)(\id \times r)=(r \times \id)(\id \times r)(r \times \id).$$
	In this case, $(X,r)$ is also called a \it{braided set}.
	A solution is called \it{non-degenerate} (\it{involutive}, or \it{square free}) if it is non-degenerate (involutive, or square free) as a quadratic set.
	If $X$ is a finite set, then $(X,r)$ is called a finite solution.

	\begin{con}
		In the followings, by a \it{solution} (of the YBE) we mean  a set-theoretic, non-degenerate and involutive solution of the YBE. We consider only \it{finite} solutions.
	\end{con}

	Each quadratic set $(X,r)$ determines a set of quadratic defining  relations $\bold{R}(r)$ defined by
	$$x y=y' x'\in \bold{R}(r) \text{ whenever } r(x,y)=(y', x').$$
	The unital semigroup $S(X,r)=\langle X;\bold{R}(r) \rangle$, with a set of generators $X$ and a set of defining relations $\bold{R}(r)$, is called the \it{semigroup associated with} $(X,r)$. The \it{group} $G(X,r)$ \it{associated with} $(X,r)$ is defined analogously. The \it{algebra associated with}  $(X, r)$  is defined as
	$A(X, r)=\kk\langle X\rangle /(\bold{R}(r))$. 
	When $(X,r)$ is a solution of the YBE, $S(X,r)$, resp. $G(X,r)$, $A(X, r)$  is called the \it{Yang-Baxter} (\it{YB}) \it{semigroup}, resp. the  \it{YB group}, the  \it{YB algebra}.

	Let $(X,r)$ be a nondegenerate quadratic set. It is well-known that $(X,r)$ is a set-theoretic solution of the YBE if and only if the following condition holds for all $x,y,z\in X$ (e.g. \cite[Lemma 2.5]{gim1}).
	$$\mathbf{l1}:  {}^{x}\left(^{y} z\right)={ }^{^{x} y}\left(^{x^{y}} z\right), \quad \mathbf{r1}: \left(x^{y}\right)^{z}=\left(x^{^yz}\right)^{y^{z}}, \quad \mathbf{lr3}: \left({ }^{x} y\right)^{\left(^{x^{y}}(z)\right)}={ }^{(x{^{^y z}})}\left(y^{z}\right).$$
	
	For every non-degenerate braided set, condition {\bf l1} implies that  the mapping $x\ra \mc{L}_x$ extends canonically to a group homomorphism 
	$$\mc{L}:G(X,r)\longrightarrow \Sym(X),$$
	where $\Sym(X)$ is symmetric group on $X$. The homomorphism $\mc{L}$ defines the \it{canonical left action} of $G(X,r)$ on the set $X$. Similarly, {\bf r1} gives the \it{canonical right action} of 
	$G(X,r)$ on the set $X$.

	The  \it{YB permutation group} $\mc{G}(X, r)$ of a finite solution $(X, r)$ is the subgroup of the symmetric 
	group $\Sym(X)$ on $X$ generated by $\{\mc{L}_x : x\in X\}$.
	
	It is known (see e.g. \cite[Lemma 1.2]{de}) that there is a one-to-one correspondence between solutions of the YBE and involutive biracks.
	
	\begin{lem}
		\begin{enumerate}
			\item If $(X,\circ, \backslash_\circ, \bullet,/_\bullet)$ is an involutive birack, then defining
			$$r(x,y)=(x\circ y, x\bullet y).$$
			We obtain that $(X,r)$ is a solution of the YB equation. 
			\item Let $(X,r)$ be a solution of the YB equation. The bijective mapping $r$ naturally defines a left
			action $\mc{L}$ and a right action $\mc{R}$ on $X$ as:
			$$r(x,y)=(\mc{L}_x(y),\mc{R}_y(x)).$$ Then  defining 
			$$\begin{array}{ll}
				x\circ y=\mc{L}_x(y),&x\backslash_\circ y =\mc{L}^{-1}_x(y),\\
				x\bullet y=\mc{R}_y(x),&x/_\bullet y=\mc{R}^{-1}_y(x).
			\end{array}$$
			We obtain an involutive birack $(  X, \circ,\backslash_{\circ}, \bullet, /_\bullet  )$.
		\end{enumerate}
	\end{lem}
	Such a correspondence allows us to treat each solution as an involutive birack.

	\section{Isotopes of quasigroups}\label{S2}
	This section serves as a preparation for Section \ref{S3}. We introduce the notion of an $\NN^p$-graded birack and  generalize the isotope construction of \cite{jpz} to this graded setting.

	We say a left quasigroup $(  X, \circ,\backslash _{\circ}  )$ is  \it{$\NN^p$-graded} if there is  a parition on $X$:
	$$X=X_{1} \cup X_{2} \cup \cdots \cup X_{p}$$
	and for any $x\in X, y\in X_i$, 	\begin{equation}\label{E7}
		x\circ y\in X_i.
	\end{equation}
	If  $x \in X_{i}$, then the degree vector is assigned as:
	$$\deg(x)=(0, \ldots, 0,1,0, \ldots, 0) \in \NN^{p},$$
	where the 1 is in the  $i$-th position. Condition (\ref{E7}) means that all left translation is degree-preserving.

	An \it{$\NN^p$-graded automorphism}  $f$ of 	$(  X, \circ,\backslash _{\circ}  )$ is an automorphism of $(X, \circ,\backslash _{\circ})$ such that $f(X_i 
	)\subseteq X_i$ for all $1\leqslant i \leqslant p$.

	Let $\phi:=\left\{\phi_{s} \mid s=1, \cdots, p\right\}$ be a sequence of commuting degree-preserving bijection of the set  $X$, i.e., for every $1\leqslant s\leqslant p$, $\phi_{s}(x)\in X_i$ for any $x\in X_i$, $1\leqslant i\leqslant p$.  For a degree vector $\alpha=\left(a_{1}, \cdots, a_{p}\right)$, set $\phi^{\alpha}=\phi_{1}^{a_{1}} \cdots \phi_{d}^{a_{p}}$. Define on the set  $X$  new binary operations:
	$$x \ast  y:=x \circ \phi^{|x|}(y)=\mc{L}_{x} \phi^{|x|}(y),$$
	$$x \backslash_{\ast} y:=\phi^{-|x|}\left(x \backslash_{\circ} y\right)=\phi^{-|x|} \mc{L}_{x}^{-1}(y).$$
	The algebra $(X,*,\backslash_*)$ is called the \it{$\phi$-isotope} of $(X,\circ, \backslash_\circ)$.
	
	It is easy to see that 
	$$x *\left(x \backslash_{*} y\right)=\mc{L}_{x} \phi^{|x|} \phi^{-|x|} \mc{L}_{x}^{-1}(y)=y$$
	and $$x \backslash_{*}(x * y)=\phi^{-|x|} \mc{L}_{x}^{-1}\mc{L}_{x} \phi^{|x|} (y)=y.$$
	So the $\phi$-isotope $(X,*,\backslash_*)$ is also a left quasigroup. Moreover, the grading of $(  X, \circ,\backslash _{\circ}  )$ induces an $\NN^p$-grading on $(X,*,\backslash_*)$.

	\begin{lem}\label{L1}
		Let $(  X, \circ, \backslash_{\circ}  )$ be a non-degenerate $\NN^p$-graded left quasigroup and let  $\phi:=\left\{\phi_{s} \mid s=1, \cdots, p\right\}$ be a sequence of commuting degree-preserving bijection of the set  $X$. Then the $\phi$-isotope $(X,*,\backslash_*)$  is also non-degenerate.
	\end{lem}
	\proof The left quasigroup $(  X, \circ, \backslash_{\circ}  )$ is $\NN^p$-graded, say \begin{equation}\label{E8}X=X_{1} \cup X_{2} \cup \cdots \cup X_{p}.\end{equation}All left translations is degree-preserving. By the non-degeneracy of the left quasigroup,  for each subset $X_s$, the mapping
	$$\begin{array}{rcl}
		X_s& \longrightarrow &X_s\\
		x &\longmapsto& \mc{L}_{x}^{-1}(x)=x \backslash_{\circ} x
	\end{array}$$
	is bijective. Each $\phi_s$ is bijective, so the mapping
	$$\begin{array}{rcl}
		X_s& \longrightarrow &X_s\\
		x &\longmapsto& \phi_s \mc{L}_{x}^{-1}(x)=x \backslash_{*} x
	\end{array}$$
	is also bijective. (\ref{E8}) is a partition of $X$. Therefore, the mapping  
	$$\begin{array}{rcl}
		T_\phi: X& \longrightarrow &X\\
		x &\longmapsto& \phi_s \mc{L}_{x}^{-1}(x)=x \backslash_{*} x
	\end{array}$$
	is bijective. The $\phi$-isotope $(X,*,\backslash_*)$  is also non-degenerate.\qed

	\begin{lem} \label{L3}Let $(X, \circ, \backslash_\circ)$ be an $\NN^p$-graded right cyclic left quasigroup, and $\phi:=\left\{\phi_{s} \mid s=1, \cdots, p\right\}$ be a sequence of commuting $\NN^p$-graded automorphisms of $(X, \circ, \backslash_\circ)$ such that 
			\begin{equation}\label{E1}
			\mc{L}_{\phi_s(x)} = \mc{L}_x.
		\end{equation}
	for every $1\leqslant s \leqslant  $ and  any $x\in X$.	Then the $\phi$-isotope $(X,*,\backslash_*)$ of $(X, \circ, \backslash_\circ)$ is also right cyclic. 
\end{lem}

\proof We need to proof that 
$$	(x\backslash_{*} y) \backslash_*(x \backslash_* z)=(y \backslash_{*} x)\backslash_{*}(y \backslash_{*} z),$$
for every $x,y,z\in X$. Note that $\left\{\phi_{s} \mid s=1, \cdots, p\right\}$ is a sequence of  automorphism of $(X, \circ, \backslash_\circ)$. Then $\mc{L}_{\phi_s(x)}=\phi_s\mc{L}_x\phi_s^{-1}$ for every $x \in X$ and $1\leqslant s \leqslant p$. In this case, Condition (\ref{E1}) is equivalent to 
\begin{equation}\label{E2}
	\phi_s\mc{L}_{x} = \mc{L}_x\phi_s.
\end{equation}
That is,  every $\phi_s$ commutes with all   left translations.

 We have the following equations hold	$$ 
\begin{array}{rl}
	&(x\backslash_{*} y) \backslash_*(x \backslash_* z)\\
	=&(\phi^{-|x|}\mc{L}^{-1}_{x} (y))\backslash_{*}(\phi^{-|x|}\mc{L}^{-1}_{x} (z))\\
	=&\phi^{-|\phi^{-|x|}(\mc{L}^{-1}_{x} (y))|}\mc{L}^{-1}_{\phi^{-|x|}\mc{L}^{-1}_{x} (y)} \phi^{-|x|}\mc{L}^{-1}_{x} (z)\\
	=&\phi^{-|y|}\mc{L}^{-1}_{\phi^{-|x|}\mc{L}^{-1}_{x} (y)} \phi^{-|x|}\mc{L}^{-1}_{x} (z)\\
	\overset{(\ref{E1})}{=}&\phi^{-|y|}\mc{L}^{-1}_{\mc{L}^{-1}_{x} (y)} \phi^{-|x|}\mc{L}^{-1}_{x} (z)\\
	\overset{(\ref{E2})}{=}&\phi^{-(|x|+|y|)}\mc{L}^{-1}_{\mc{L}^{-1}_{x} (y)} \mc{L}^{-1}_{x} (z)\\
	=&\phi^{-(|x|+|y|)}\mc{L}^{-1}_{\mc{L}^{-1}_{y} (x)} \mc{L}^{-1}_{y} (z)\;\;\tiny{\text{($(X, \circ, \backslash_\circ)$ is right cyclic)}}\\
	\overset{(\ref{E2})}{=}&\phi^{-|x|}\mc{L}^{-1}_{\mc{L}^{-1}_{y} (x)} \phi^{-|y|}\mc{L}^{-1}_{y} (z)\\
	\overset{(\ref{E1})}{=}&\phi^{-|x|}\mc{L}^{-1}_{\phi^{-|y|}\mc{L}^{-1}_{y} (x)} (x)\phi^{-|y|}\mc{L}^{-1}_{y} (z)\\
	=&\phi^{-|\phi^{-|y|}\mc{L}^{-1}_{y} (x)|}\mc{L}^{-1}_{\phi^{-|y|}\mc{L}^{-1}_{y} (x)} (x)\phi^{-|y|}\mc{L}^{-1}_{y} (z)\\
	=&(\phi^{-|y|}\mc{L}^{-1}_{y} (x))\backslash_{*}(\phi^{-|y|}\mc{L}^{-1}_{y} (z))\\
	=&(y \backslash_{*} x)\backslash_{*}(y \backslash_{*} z).
\end{array}		$$
The third and the last third equation follow from the fact that $(X, \circ, \backslash_\circ)$ is  an $\NN^p$-graded  left quasigroup and $\phi:=\left\{\phi_{s} \mid s=1, \cdots, p\right\}$ is a sequence of $\NN^p$-graded automorphism of $(X, \circ, \backslash_\circ)$.\qed

We know that involutive biracks are in one-to-one correspondence with non-degenerate right cyclic left quasigroups. Now we are ready to  define the isotope of an $\NN^p$-birack. 
	
	A birack $(X,\circ,\backslash_\circ,\bullet, \backslash_\bullet)$ is called  \it{$\NN^p$-graded} if there is  a parition on $X$:
	$$X=X_{1} \cup X_{2} \cup \cdots \cup X_{p}$$ 
	and for any $x\in X, y\in X_i$, 	
	$$x\circ y\in X_i,\text{ and } y\bullet x \in X_i.$$
	
\begin{rem}
If $(X,\circ,\backslash_\circ,\bullet, /_\bullet)$ is an involutive birack, then $(X,\circ,\backslash_\circ)$ is a non-degenerate right cyclic left quasigroup.  The birack $(X,\circ,\backslash_\circ,\bullet, /_\bullet)$ is an $\NN^p$-graded birack if and only if the $(X,\circ,\backslash_\circ)$ is an $\NN^p$-graded left quasigroup. 
\end{rem}	
	
	Let  $\left(X, \circ, \backslash_{\circ}, \bullet, /_\bullet\right)$  be an $\NN^p$-graded involutive birack. Then $(  X, \circ, \backslash_{\circ}  )$ is  an $\NN^p$-graded non-degenerate right cyclic left quasigroup. Let  $\phi:=\left\{\phi_{s} \mid s=1, \cdots, p\right\}$  be   a sequence of commuting degree-preserving bijection of the  set  $X$  such that the  $\phi$-isotope  $\left(X, *, \backslash_{*}\right)$  of $(  X, \circ \backslash_{\circ}  )$ is right cyclic. Then  the left quasigroup  $\left(X, *, \backslash_{*}\right)$  uniquely determines the involutive birack $ \left(X, *, \backslash_{*}, \diamond, /_\diamond\right) $ as follows
	\begin{equation}
		x * y=\mc{L}_{x} \phi^{|x|}(y),
	\end{equation}
	\begin{equation}	x \diamond y=(x * y) \backslash_{*} x=\phi^{-|y|} \mc{L}_{\mc{L}_{x} \phi^{|x|}(y)}^{-1}(x)=\phi^{-|y|} (x\bullet \phi^{|x|}(y)).\end{equation}
	The birack $ \left(X, *, \backslash_{*}, \diamond, /_\diamond\right) $ is called the \it{$\phi$-isotope} of $\left(X, \circ, \backslash_{\circ}, \bullet, /_\bullet\right)$. It is easy to see that $ \left(X, *, \backslash_{*}, \diamond, /_\diamond\right) $ is also an $\NN^p$-graded birack. 
	
	\begin{lem}\label{L2}
		Let $(X, \circ, \backslash_\circ, \bullet,/_\bullet)$ be an $\NN^p$-graded involutive birack satisfying Condition {\bf lri}, and $\phi:=\left\{\phi_{s} \mid s=1, \cdots, p\right\}$  be   a sequence of commuting degree-preserving bijection of the  set  $X$. Consider the $\phi$-isotope $(X, *, \backslash_{*}, \diamond, /_\diamond)$ of the set $X$. Then $(X, *, \backslash_{*}, \diamond, /_\diamond)$ satisfies Condition {\bf lri} if and only if  $(X, \circ, \backslash_\circ, \bullet,/_\bullet)$ satisfies Condition (\ref{E1}).
			In this case, 	$$y\diamond x =x\backslash_*y=\phi^{-|x|} \mc{L}_{x}^{-1}(y),$$ for any $x,y \in X$.
	\end{lem} 
	\proof It is clear that if $(X, *, \backslash_{*}, \diamond, /_\diamond)$ satisfies Condition {\bf lri}, then 
	$$y\diamond x =x\backslash_*y=\phi^{-|x|} \mc{L}_{x}^{-1}(y),$$ for any $x,y \in X$. 
	
	 Indeed, for every $1\leqslant s\leqslant p$, $y\in X_s$ and $x \in X$, we have
	$$\begin{array}{rcl}
		y\diamond x =x\backslash_*y	&\iff&\phi^{-|x|} \mc{L}_{\mc{L}_{y} \phi^{|y|}(x)}^{-1}(y)=\phi^{-|x|} \mc{L}_{x}^{-1}(y)\\
		&\iff&\phi^{-|x|}(y \bullet \phi^{|y|}(x))=\phi^{-|x|} (y\bullet x)\\
		&\iff& y \bullet \phi^{|y|}(x)= y\bullet x\\
		&\iff& \mc{R}_{\phi^{|y|}(x)} = \mc{R}_x\\
		&\iff& \mc{L}_{\phi_s(x)} = \mc{L}_x.
	\end{array}$$
Therefore, $(X, *, \backslash_{*}, \diamond, /_\diamond)$ satisfies Condition {\bf lri} if and only if  $(X, \circ, \backslash_\circ, \bullet,/_\bullet)$ satisfies Condition (\ref{E1})
\qed

	\begin{prop}\label{P1}
		Let $(X,\circ, \backslash_\circ, \bullet,/_\bullet)$  be an $\NN^p$-graded involutive birack satifying Condition {\bf lri}, and $\phi:=\left\{\phi_{s} \mid s=1, \cdots, p\right\}$  be a sequence of commuting graded automorphisms of $(X, \circ, \backslash_\circ, \bullet,/_\bullet)$ satisfying Condition (\ref{E1}).  Then the  $\phi$-isotope $(X, *, \backslash_{*}, \diamond, /_\diamond)$ of $(X,\circ, \backslash_\circ, \bullet,/_\bullet)$ is also an  $\NN^p$-graded involutive birack satifying Condition {\bf lri}.
	\end{prop}
	\proof $(X,\circ, \backslash_\circ, \bullet,/_\bullet)$ is an $\NN^p$-graded involutive birack. Then $(X,\circ, \backslash_\circ)$ is a non-degenerate right cyclic left quasigroup by Lemma \ref{L4}. By Lemma \ref{L1} and \ref{L3}, the $\phi$-isotope $(X,*,\backslash_*)$ is also non-degenerate right cyclic. It determines the $\NN^p$-graded involutive birack $(X, *, \backslash_{*}, \diamond, /_\diamond)$. By Lemma \ref{L2},  $(X, *, \backslash_{*}, \diamond, /_\diamond)$ also satisfies Condition {\bf lri}. \qed


	\section{Isotopes and Zhang twists}\label{S3}
In this seciton, we will show the relation between the isotopes of involutive biracks and the Zhang twists of algebras. 

	Let $(X,\circ, \backslash_\circ, \bullet,/_\bullet)$ be an involutive birack. It canonically determines an algebra $A=A(X,\bf{R^\circ})$ generated by $X$ with quadratic defining relations $\bf{R}^\circ$ defined by 
	$$xy=(x\circ y)(x\bullet y) $$ 
	for any $x,y\in X$. 
	
	If $(X,\circ, \backslash_\circ, \bullet,/_\bullet)$ is an $\NN^p$-graded birack, then the grading on $X$ naturally defines a grading on $A$.  The algebra $A$ is an  $\NN^p$-graded algebra. 
	
	\begin{rem}
		For an ($\NN^p$-graded)  involutive birack $(X,\circ, \backslash_\circ, \bullet,/_\bullet)$, the ($\NN^p$-graded)  algebra  $A=A(X,\bf{R^\circ})$  is just the YB algebra  of the solution corresponding to the  involutive birack.
	\end{rem}
	Let  $B$  be any  $\mathbb{N}^{p}$-graded algebra. Let  $\phi:=\left\{\phi_{s} \mid s=1, \cdots, p\right\}$  be a sequence of commuting  $\mathbb{N}^{p}$-graded algebra automorphisms of  $B$. Then the \it{Zhang twist}  $B^{\phi}$ of $B$ by $\phi$ is definied as follows:  $B^{\phi}=B$  as a  $\mathbb{N}^{p}$-graded vector space, and the new multiplication $\star$ of  $B^{\phi}$  is determined by
	$$a \star b=a \phi^{|a|}(b)$$
	where  $\phi^{|a|}=\phi_{1}^{a_{1}} \cdots \phi_{p}^{a_{p}}$  if the degree of  $a$  is  $|a|=\left(a_{1}, \cdots, a_{p}\right)$.

	Let $(X,\circ, \backslash_\circ, \bullet,/_\bullet)$  be an $\NN^p$-graded involutive birack satifying condition {\bf lri}. It determines an $\mathbb{N}^{p}$-graded  algebra $A=A(X,\bf{R}^\circ)$.  Let $\phi:=\left\{\phi_{s} \mid s=1, \cdots, p\right\}$  be a sequence of commuting graded automorphism of $(X, \circ, \backslash_\circ, \bullet,/_\bullet)$ satisfying Condition (\ref{E1}). It is easy to see that $\left\{\phi_{s} \mid s=1, \cdots, p\right\}$ naturally induce a sequence of commuting graded automorphism on the algebra $A$.  By Proposition \ref{P1}, the $\phi$-isotope  $(X,*, \backslash_*, \diamond,/_\diamond)$  of $(X,\circ, \backslash_\circ, \bullet,/_\bullet)$ is also  an $\NN^p$-graded involutive birack satifying condition {\bf lri}. It determines an $\NN^p$-graded algebra $A^*=A^*(X,\bf{R^*})$.

	\begin{thm}\label{T1} Keep the notations as above,  the algebra $A^*=(A^*,X,\bf{R}^*)$ obtained by the isotope $(X,*, \backslash_*, \diamond,/_\diamond)$ is isomorphic to the Zhang twist $A^\phi$ of $A$ by $\phi$.
	\end{thm}
	\proof  For the algebra $A^*=(A^*,X,\bf{R}^*)$, the relation $\bf{R}^*$ is defined by  
	$$\begin{array}{rcl}
		xy&=&(x * y)(x\diamond y)=(x * y)(y \backslash_{*}x ) \\
		&=&	\mc{L}_x(\phi^{|x|}(y)) \phi^{-|y|}(\mc{L}^{-1}_y(x)).
	\end{array}$$
	We need to prove that $A^\phi$ is isomorphic to the algebra generated by $X$ with relation ${\bf R^*}$. We have
	$$\begin{array}{rcl}
		x\star y&=&x\phi^{|x|}(y)\\
		&=&(x\circ \phi^{|x|}(y)) (\phi^{|x|}(y)\backslash_{\circ} x)\\
		&=&\mc{L}_x(\phi^{|x|}(y))\mc{L}^{-1}_{\phi^{|x|}(y)}(x)\\
		&\overset{({\bf lri})}{=}&\mc{L}_x(\phi^{|x|}(y))\mc{L}^{-1}_{y}(x)\\
		&=&\mc{L}_x(\phi^{|x|}(y))\star \phi^{-|\mc{L}_x(\phi^{|x|}(y))|}(\mc{L}_y(x))\\
		&=&\mc{L}_x(\phi^{|x|}(y))\star \phi^{-|y|}(\mc{R}_y(x)).
	\end{array}
	$$
The last equation follows from that each $\phi_s$ ($1\leqslant s \leqslant p$) and all left translation preserve degree. Therefore, the Zhang twist of $A^\phi$ is isomorphic to the algebra $A^*$.\qed

Since the Zhang twist preserves many important properties of graded algebras, with the help of Theorem \ref{T1}, we can reduce  the study of a structurally complicated YB algebra to that of a simpler one via isotopy. The following is an example,  we will show that YB algebras of distributive solutions are Zhang twists of polynomial algebras. 

	\begin{eg} 
		
		A birack $(X,\circ, \backslash_\circ, \bullet,/_\bullet)$ is called \it {distributive}, if  for every $x, y, z \in X$,
		$$x \circ(y \circ z)=(x \circ y) \circ(x \circ z),$$
		$$(y \bullet z) \bullet x = (y \bullet x) \bullet (z \bullet x).$$
		
		If the birack $(X,\circ, \backslash_\circ, \bullet,/_\bullet)$  involutive, then it is distributive if and only if  for every $x, y, z \in X$, \begin{equation}\label{E11}
			x \circ(y \circ z)=(x \circ z) \circ(x \circ z),
		\end{equation} i.e., 	 $\mc{L}_x\mc{L}_y=\mc{L}_{x\circ y}\mc{L}_x$ \cite[Corollary 5.7]{jpz}. Following from \cite[Corollary 5.4]{jpz} and \cite[Lemma 7.1]{gi3}, we know that  an involutive distributive birack satisfies Condition {\bf lri}.
		
		A solution of YBE is called \it{distributive}, if it corresponds to a distributive  involutive birack.
		
		Let $(X,\circ, \backslash_\circ, \bullet,/_\bullet)$ be an involutive  distributive birack. There is an equivalence $\sim$ defined on $X$:
		$$x\sim y  \quad\text{ if and only if } \quad\mc{L}_x=\mc{L}_y.$$
		This equivalence gives a partition on $X$. That is, 
		\begin{equation}\label{E9}
			X=\bigcup_{x\in X}[x],
		\end{equation}
		where $[x]=\{a\in X|x\sim a\}.$
		By Lemma 3.1 and Corollary 5.4 in \cite{jpz},  an involutive  distributive birack is  \it{2-redective}, i.e., 
		\begin{equation}\label{E10}
			(x\circ y)\circ z=  y\circ z\;\;(\iff \mc{L}_{\mc{L}_x(y) }=\mc{L}_y)
		\end{equation}
		for any $x,y,z\in X$. Then for any $y\in [x]$ and $x\in X$, $\mc{L}_x(y)\in [x]$, (\ref{E9}) gives a grading on the birack.  Assume that $X$  is partitioned into $p$  equivalence classes, numbering them as $X_1,\cdots, X_p$. The birack $(X,\circ, \backslash_\circ, \bullet,/_\bullet)$ is an $\NN^p$-graded birack.  
		
		In each equivalence class $X_s$, select a representative element $x_s$.  We have $\mc{L}_y=\mc{L}_{x_s}$, for  any $y\in X_s$. We denote the inverse of these equal mappings as $\phi_s$, i.e. $\phi_s=\mc{L}_{x_s}^{-1}$.  Directly from (\ref{E11}) and (\ref{E10}), we obtain that all left translations are automorphisms of the birack satisfying Condition (\ref{E1}). Moreover, the permutation group generated by the left transformations are abelian. So $\phi:=\left\{\phi_{s} \mid s=1, \cdots, p\right\}$  is a sequence of commuting graded automorphisms of the birack satisfying Condition (\ref{E1}). The 
		$\phi$-isotope of the birack $(X,\circ, \backslash_\circ, \bullet,/_\bullet)$ is just the projection birack, i.e. the birack whose all left and right translations are  the identity. Since involutive distributive biracks satisfy Condition {\bf lri}, we only need to check for the left translations. Indeed, 
		for any $x,y\in X$, 
		$$	x * y=\mc{L}_{x} \phi^{|x|}(y)=\mc{L}_{x} \mc{L}_{x}^{-1}(y)=y.$$
		The projection birack corresponds the trivial solution of the YBE, whose YB algebra is just the polynomial algebra. Therefore, we obtain the following proposition. 
		\begin{prop}
			The YB algebra of  a distributive solution $(X,r)$ of the YBE is a Zhang twist of the polynomial algebra over the set $X$. 
		\end{prop}

	\end{eg}

	\subsection*{Acknowledgement}  The  authors are supported by a grant from NSFC (No. 12371017).

	\vspace{5mm}
	
	\bibliography{}

\end{document}